\documentclass[12pt]{amsart}
\usepackage{palatino,hhline, bbm}
\usepackage{arydshln} 
\usepackage{amsthm, amsmath}
\usepackage{amssymb} 
\usepackage{amscd}
\usepackage{mathrsfs}
\renewcommand{\mathcal}{\mathscr}
\usepackage[hmargin=2.8cm, vmargin=2.8cm]{geometry}
\usepackage[breaklinks,bookmarksopen,bookmarksnumbered]{hyperref}
\usepackage[all]{xy}

\theoremstyle{plain}

\newtheorem{lem}{Lemma}
\newtheorem{prop}{Proposition}
\newtheorem*{cor}{Corollary}

\theoremstyle{remark}

\newtheorem{conj}{Conjecture}

\newcommand\pr{\noindent\textit{Proof}.\hskip2mm}

\newcommand\rond{\kern 1pt{\scriptstyle\circ}\kern 1pt}

\def\lr#1{\langle {#1} \rangle}
\newcommand\End{\operatorname{End}}

\newcommand\codim{\operatorname{codim}}
\newcommand\Hom{\operatorname{Hom}}
\newcommand\im{\operatorname{Im}}
\newcommand\Ker{\operatorname{Ker}}

\newcommand\Tr{\operatorname{Tr}}
\newcommand\pp{\scriptscriptstyle\bullet}
\newcommand{\mo}{\smallsetminus}
\newcommand\rk{\operatorname{rk} }

\newcommand\Q{\mathbb{Q}}

\newcommand\C{\mathbb{C}}
\renewcommand\P{\mathbb{P}}
\newcommand\N{\mathbb{N}}

\newcommand\G{\mathbb{G}}
\renewcommand\O{\mathcal{O}}

\renewcommand\ss{\mathsf{S}}

\newcommand\iso{\vbox{\hbox to .8cm{\hfill{$\scriptstyle\sim$}\hfill}
\nointerlineskip\hbox to .8cm{{\hfill$\longrightarrow $\hfill}} }}
\newcommand\bir{\vbox{\hbox to .8cm{\hfill{$\scriptstyle\sim$}\hfill}
\nointerlineskip\hbox to .8cm{{\hfill$\dasharrow $\hfill}} }}
\newcommand\abs[1]{\lvert {#1}\rvert}

\begin{document}
\title{The algebra of symmetric tensors on smooth projective varieties}
\author[Arnaud Beauville]{Arnaud Beauville}
\address{Universit\'e C\^ote d'Azur\\
CNRS -- Laboratoire J.-A. Dieudonn\'e\\
%UMR 7351 du CNRS\\
Parc Valrose\\
F-06108 Nice cedex 2, France}
\email{arnaud.beauville@unice.fr}
 \author{Jie Liu}
\address{Institute of Mathematics\\
	Academy of Mathematics and Systems Science\\
	Chinese Academy of Sciences\\
	Beijing\\
	100190\\
	China}
\email{jliu@amss.ac.cn}

\begin{abstract}
We discuss in this  note  the $\C$-algebra $H^0(X,\ss^{\pp}T_X)$ for a smooth complex projective variety $X$. We compute it in some simple examples, and give a sharp bound on its Krull dimension. Then we propose a conjectural characterization of non-uniruled projective manifolds with pseudo-effective tangent bundle. \end{abstract}
\thanks{A.~Beauville is indebted to Feng Shao for several useful comments and references. J.~Liu would like to thank S.~Druel for useful communication. J.~Liu is supported by the National Key Research and Development Program of China (No. 2021YFA1002300), the NSFC grant (No. 12288201) and the CAS Project for Young Scientists in Basic Research (No. YSBR-033)}
\maketitle 
\section{Introduction}
Let $X$ be a smooth complex projective variety. We will denote by $T^*X$ its cotangent bundle, and by $\P T^*X$ its projectivization. In this note we are interested in the graded $\C$-algebra 
\[S(X):=\bigoplus_{p\geq 0}H^0(X,\ss^pT_X)= \O(T^*X)=\bigoplus_{p\geq 0} H^0(\P T^*X,\O_{\P T^*X}(p))\,.\]
($\O(T^*X)$ is the algebra of regular functions on $T^*X$, with the grading defined by the linear action of $\C^*$.)

Despite its simple definition, this is an intriguing object, which is usually quite complicated, even for a variety as simple as the quadric (see Proposition \ref{quadric} below). While the algebra $\bigoplus_{p\geq 0}H^0(X,\ss^p\Omega ^1_X)$ has been extensively studied, starting with Sakai's work \cite{S}, this is not the case of $S(X)$.  We will describe it in some particular cases in \S 2 and 3. Then we give a sharp bound on the Krull dimension of $S(X)$ (\S 4). Finally we propose a conjectural characterization of non-uniruled projective manifolds with pseudo-effective tangent bundle, which holds in dimension $\leq 5$ (\S 5).
\medskip	
 
 We would like to dedicate this note to the memory of   Gang Xiao. Gang was (informally) a student of the first author  in Orsay at the beginning of the 80's, then later his colleague and friend in Nice, till his untimely death in 2014. 
  
 \smallskip	
 \noindent\textbf{Notations}: We work over the complex numbers. If $X$ is a variety endowed with an action of $\C^*$, we denote by $\O(X)$ the $\C$-algebra of regular functions on $X$, with the grading defined by the $\C^*$-action.  By a vector space we mean a complex, finite-dimensional vector space.
 
\smallskip	 
\section{Some examples}
\subsection{Abelian varieties}\label{ab}
We start with a trivial case: if $X$ is an abelian variety of dimension $n$, we have $T_X\cong \O_X^n$, hence $S(X)$ is a polynomial algebra in $n$ variables.

\smallskip	
\subsection{Projective space}
Let $V$ be a   vector space. We let $I\in V\otimes V^*$ be the image of the identity by the isomorphism $\End(V)\iso V\otimes V^*$.
\begin{prop}\label{proj}
The graded algebra $S(\P(V))$ is isomorphic to the quotient of $\displaystyle\bigoplus_{d\geq 0}\bigl(\ss^dV\otimes \ss^dV^*\bigr)$ by the ideal generated by $I$. 
\end{prop}
\pr The projective cotangent bundle $\P T^*_{\P(V)}$ can be identified with the incidence hypersurface $Z\subset \P(V)\times \P(V^*)$ consisting of pairs $(x,H)$ with $x\in H$; the tautological line bundle $\O_{Z}(1)$ is induced by $\O_{\P(V)}(1)\boxtimes \O_{\P(V^*)}(1)$. The Proposition follows from the exact sequence
\[0\rightarrow \O_{\P(V)}(d-1)\boxtimes \O_{\P(V^*)}(d-1)\xrightarrow{\ \times I\ }\O_{\P(V)}(d)\boxtimes \O_{\P(V^*)}(d)\rightarrow \O_Z(d)\rightarrow 0\,.\qed\]

\subsection{Rational homogeneous manifolds}

In this section we will use some general facts about nilpotent orbits, which can be found for example in \cite{Fu}.

Let $X=G/P$, where $G$ is a reductive algebraic group and $P$ a parabolic subgroup. We denote by $\mathfrak{g}$ and $\mathfrak{p}$ their Lie algebras, and by $\mathfrak{n}$ the nilradical of $\mathfrak{p}$. The Killing form of $\mathfrak{g}$ provides an isomorphism of $G$-modules $\mathfrak{n}\iso(\mathfrak{g}/\mathfrak{p})^*$; using this we identify
the cotangent bundle $T^*(G/P)$ to the homogeneous bundle $G\times ^P\mathfrak{n}$. Associating to a pair $(g,N)$ in $G\times \mathfrak{n}$ the element $\operatorname{Ad}(g)\cdot N $ of $\mathfrak{g}$ defines a generically finite, $\C^*$-equivariant map $\pi :T^*(G/P)\rightarrow \mathfrak{g}$, whose image $ \mathscr{N} $ is the closure of a nilpotent orbit. 

We will consider the case where the induced map  
 $\bar{\pi }:T^*(G/P)\rightarrow \mathscr{N} $ is birational. In this case $\bar{\pi }$ is a resolution of the normalization $\tilde{\mathscr{N}}$ of $\mathscr{N}$, and we have $S(X)=\O(\tilde{\mathscr{N} })$. For $G=\operatorname{GL}(n) $ all parabolic subgroups have this property, and $\mathscr{N} $ is normal, so $S(X)=\O( \mathscr{N} )$. In the other classical cases there is a precise description of the 
 parabolic subgroups for which $\bar{\pi} $ is birational
 \cite[3.3]{Fu}; we will content ourselves with the example of quadrics.

\smallskip	
\subsection{Flag varieties, Grassmannians}

Let $V$ be a  vector space, and let $(0)=V_0\subset V_1\subset\allowbreak\ldots \subset V_{s+1}=V$ be a (partial) flag in $V$. The stabilizer $P$ of this flag is a parabolic subgroup of $\operatorname{GL}(V) $, and all parabolics are  obtained in this way. The variety $G/P$ is the variety of flags $(0)=F_0\subset F_1\subset\ldots \subset F_{s+1}=V$ with $\dim F_i=\dim V_{i}$.

The Lie algebra $\mathfrak{p}$ is  the stabilizer of $(V_i)$ in $\End(V)$, and its nilradical $\mathfrak{n}$ is the subspace of $u\in \End(V)$ satisfying $u(V_{i+1})\subset V_i$ for $0\leq i\leq s$. Therefore $\mathscr{N}$ is the subvariety of endomorphisms $u\in\End(V)$ for which there exists a flag $(F_i)$ in $G/P$ with $u(F_{i+1})\subset F_i$ for $0\leq i\leq s$. 

Let us spell out this in the case of the Grassmannian $\G:=\G(r,V)$ of $r$-dimensional subspaces of $V$. We put $n:=\dim V$.
\begin{prop}
$S(\G(r,V))=\O(\mathscr{N})$, where $\mathscr{N}\subset \End(V)$ is the subvariety of endomorphisms $u$ satisfying $u^2=0$ and $\operatorname{rk}  u\leq \min\{r,n-r\} $.
\end{prop}
\pr Since $\G(r,V)\cong \G(n-r,V)$, we can assume $r\leq n/2$. By the previous discussion, $\mathscr{N}$ consists of endomorphisms $u$ for which there exists an $r$-dimensional subspace $W\subset V$ with $u(V)\subset W$ and $u(W)=0$, that is, $\im u\subset W\subset\Ker u$. This implies $u^2=0$ and $\operatorname{rk}u\leq r $; conversely, if this is satisfied, we have $\im u\subset \Ker u$ and $\dim\Ker u=n-\operatorname{rk}u\geq n-r\geq r $, so any  $r$-dimensional subspace $W$ with $\im u\subset W\subset\Ker u$ does the job.\qed

\smallskip	
\noindent\emph{Remarks}$.-$ 1) Taking $r=1$ we recover Proposition \ref{proj}.

\noindent 2) If $r=\lfloor\frac{n}{2} \rfloor$ the condition $u^2=0$  implies  $\operatorname{rk}u \leq r$, so $\mathscr{N}$ is simply the variety of square zero endomorphisms of $V$.

\smallskip	
\subsection{Quadrics}
Let $V$ be a  vector space, and let $q$ be a non-degenerate quadratic form on $V$, defining a quadric $Q:=V(q)$ in $\P(V)$. 

\begin{prop}\label{quadric}
$S(Q)$ is isomorphic to the quotient of the homogeneous coordinate ring of $\G(2,V)\subset \P(\bigwedge^2V)$ by the ideal generated by $\wedge^2q$.

\end{prop}
\pr Let $\ell$ be an isotropic line in $V$ and let $P$ be the stabilizer of $\ell$, so that $Q= \operatorname{O}(V)/P $. The Lie algebra $\mathfrak{o}(V)$ consists of endomorphisms of $V$ which are skew-symmetric (with respect to $q$), and $\mathfrak{p}$ is the stabilizer of $\ell$ in $\mathfrak{o}(V)$.

The nilradical $\mathfrak{n}$ of $\mathfrak{p}$ consists of skew-symmetric endomorphisms $u$ such that $u(\ell^{\perp})\subset \ell$ and $u(\ell)=0$. Such a map is of the form 
\begin{equation}\label{map}
x\,\mapsto\, q(w,x)v-q(v,x)w, \mbox{ where } v\in\ell \mbox{ and } w\in \ell^{\perp} \,.
\end{equation}
Varying $\ell$, we see that $\mathscr{N}$ consists of the maps of the form (\ref{map}) such that the restriction of $q$ to $\lr{v,w}$ has rank $\leq 1$. Such  maps correspond bijectively to decomposable bivectors $v\wedge w\in \bigwedge^2V$, and the condition on $q$ can be written $\wedge^2q(v\wedge w)=0$. This implies the Proposition.\qed

\smallskip	
\subsection{Intersection of two quadrics}
The following result is proved in \cite{BEHLV}:
\begin{prop}\label{2Q}
Let $X\subset \P^{n+2}$ be a smooth complete intersection of two quadrics, with $n\geq 2$. Then 
$S(X)$ is a polynomial algebra in $n$ variables of degree $2$.
\end{prop}
It is somewhat surprising that the answer is much simpler in this case that for a single quadric.

\smallskip	
\subsection{Completely integrable systems}
Let $V$ be a graded vector space, endowed with  the associated $\C^*$-action. Suppose that we have a $\C^*$-equivariant morphism \break$\Phi :T^*X\rightarrow V$ whose general fiber is of the form $Y\mo Z$, where $Y$ is a complete variety and $Z$ a closed subvariety of codimension $\geq 2$. Then the functions on $T^*X$ are constant on the fibers of $\Phi $, hence the homomorphism 
$\Phi ^*: \O(V)=\ss^{\pp}V^*\rightarrow \O(T^*X)=S(X)$ is an isomorphism of graded algebras. 

A famous example of this situation is given by the \emph{Hitchin fibration} \cite{Hi}. 
Let $C$ be a curve of genus $g\geq 2$. We fix coprime integers $r,d\geq 1$, and consider  the  moduli space $\mathscr{M}$ of stable vector bundles on $C$ of rank $r$ and degree $d$.
It is a smooth projective variety. By deformation theory the tangent space $T_E(\mathscr{M})$ at a point $E$ of $\mathscr{M}$ identifies with $H^1(C,\mathscr{E}nd(E))$; by Serre duality, its dual $T^*_E\mathscr{M}$ identifies with $\Hom(E,E\otimes K_C)$. Let $ V$ be the graded vector space $\bigoplus_{i=1}^rH^0(C, K_C^{i})$ (with $\deg H^0(C, K_C^{i})=i$). For $u\in \Hom(E,E\otimes K_C)$, we have $\Tr\wedge^i u\in H^0(C,K^i_C) $. 
Associating to $u$ the  vector $\Tr u+\ldots +\Tr \wedge^r u$ gives a $\C^*$-equivariant map $\Phi :T^*\mathscr{M}\rightarrow V$.
\begin{prop}
The homomorphism $\Phi^*:\O(V)=\ss^{\pp}V^*\rightarrow \O(T^*\mathscr{M})=S(\mathscr{M}) $ is an isomorphism.
\end{prop}
\pr $T^*\mathscr{M}$ admits an open embedding into the moduli space $\mathscr{H} $of stable Higgs bundles (of rank $r$ and degree $d$), and $\Phi $ extends to a proper map $\bar{\Phi }:\mathscr{H}\rightarrow V$ \cite{Hi}.The codimension of $\mathscr{H}\mo T^*\mathscr{M}$ is $\geq 2$ \cite[Theorem II.6]{Fa},  hence  $\codim \bar{\Phi }^{-1}(v)\mo \Phi ^{-1}(v)\geq 2$ for $v$ general in $V$. By the previous remarks this implies the result.\qed

\medskip	
There are a number of variations on this theme. First of all, one can fix a line bundle $L$ of degree $d$ on $X$ and consider the subspace $\mathscr{M}_L$ of $\mathscr{M}$ parameterizing the vector bundles  $E$ with $\det E=L$; then $\Phi $ maps $\mathscr{M}_L$ onto the graded subspace $V_0:= \bigoplus_{i=2}^rH^0(C, K_C^{i})$ of $V$, and we get as before an isomorphism of $S(\mathscr{M}_L)$ with $\ss^{\pp}V_0^*$. Note that in the case $g=r=2$ $\mathscr{M}_L$ is a complete intersection of two quadrics in $\P^5$, so we recover the  case $n=3$ of Proposition \ref{2Q}.

We can also consider the moduli space $\mathscr{M}_{\operatorname{par} }$ of  stable parabolic vector bundles on $C$ of rank $r$, degree $d$ and weights $\alpha $, with a parabolic structure along a divisor $D=p_1+\ldots +p_s$ --- we refer for instance to \cite{BGL} for the precise definitions. For generic weights $\mathscr{M}_{\operatorname{par} }$ is smooth and projective; the Hitchin map $\Phi : T^*\mathscr{M}_{\operatorname{par} }\rightarrow V_{\operatorname{par} }$ takes its values in the vector space $ V_{\operatorname{par} }:=\bigoplus_{i=1}^r H^0(C,K_C((i-1)D))$. It extends to a proper map from the moduli space $\mathscr{H}_{\operatorname{par} }$ of parabolic Higgs bundle to $V_{\operatorname{par} }$, and $\mathscr{M}_{\operatorname{par} }\mo T^*\mathscr{M}_{\operatorname{par} }$ has codimension $\geq 2$ provided $g\geq 4$, or $g=3$ and $r\geq 3$, or $g=2$ and $r\geq 5$ \cite[Proposition 5.10]{BGL}. If this holds, we get as before an isomorphism $\ss^{\pp}V_{\operatorname{par} }^*\iso S(\mathscr{M}_{\operatorname{par} })$.

 \smallskip	
\subsection{An example: ruled surfaces}
Contrary to what the previous examples might suggest, $S(X)$ is \emph{not} invariant under deformation of $X$; a typical example is provided by  ruled surfaces.
Let $C$ be a curve of genus $\geq 2$, and $E$ a stable rank 2 vector bundle on $C$ with trivial determinant\footnote{Such a bundle is isomorphic to its dual, so we will not bother to distinguish them.}. We put  $X=\P_C(E)$.
\begin{prop}
For general $E$ we have $S(X)=\C$.
\end{prop}
\pr   Denote by $p:X\rightarrow C$ the structure map and by $\O_X(1)$ the tautological line bundle. The exact sequence
\[0\rightarrow \O_X(2)\rightarrow T_X\rightarrow p^*T_C\rightarrow 0\,.\]
gives rise to exact sequences
\begin{equation}\label{ex}
0\rightarrow \O_X(2p)\rightarrow \ss^pT_X\rightarrow \ss^{p-1}T_X\otimes p^*T_C \rightarrow 0\,.\end{equation}
We claim that  $H^0(X,\ss^{p-1}T_X\otimes p^*T_C)=0$. Indeed we get from (\ref{ex})  exact sequences
\[0\rightarrow \O_X(2q)\otimes p^*T_C^{r}\rightarrow \ss^qT_X \otimes p^*T_C^{r}\rightarrow \ss^{q-1}T_X\otimes  p^*T_C^{r+1} \rightarrow 0\,.\]
We have $H^0(X,\O_X(2q)\otimes p^*T_C^{r})=H^0(C, \ss^{2q}E\otimes T_C^{r}))=0$ for $r\geq 1$, because $\ss^{2q}E$ is semi-stable \cite[ch.\,I, Theorem 10.5]{Ha} and $\deg T_C<0$. Since $H^0(C,   T_C^{q+1}))=0$, we get by induction $H^0(X,\ss^{q}T_X\otimes p^*T_C)=0 $, hence  (\ref{ex}) gives  isomorphisms \begin{equation}\label{iso}
H^0(X,\ss^p T_X)\cong H^0(X,\O_X(2p))\cong H^0(C,\ss^{2p}E)\,.
\end{equation} 
Now for general $E$ the bundles $\ss^qE$ are stable \cite[\emph{loc.\ cit.}]{Ha}, so $H^0(X, \ss^pT_X)=0$ for $p>0$.\qed

\medskip		
For special bundles $E$ the algebra $S(X)$ can be quite nontrivial. If $E$ is \emph{unstable} the tangent bundle $T_{X}$ is big \cite{Ki}, hence $S(X)$ has Krull dimension $4$. This does not hold if $E$ is stable, but one can get interesting algebras of dimension 2. Let $V$ be a 2-dimensional Hermitian space, and let $G $ be a finite subgroup of $\operatorname{SU}(V) $, acting irreducibly on $V$. Recall that $G$ is the pull-back by the covering map $\operatorname{SU}(2)\rightarrow \operatorname{SO}(3)  $ of 
a group $\bar{G}$ isomorphic to the dihedral group $D_n$ or to $\mathfrak{A}_4,\mathfrak{S}_4$ or $\mathfrak{A}_5$.

Given an \'etale Galois covering $\pi :\tilde{C}\rightarrow C $ with group $G$,   the vector bundle \allowbreak$E_{\pi }:= \tilde{C}\times ^{G}V $ on $C$ is stable, of rank 2, with trivial determinant.
The space $H^0(C,\ss^pE_{\pi })$ is canonically isomorphic to the $G$-invariant subspace of $\ss^pV$. Note that this is zero if $p$ is odd, since $G$ contains the element $-1_V$. Therefore it follows from (\ref{iso}) that $S(X)$ \emph{is isomorphic to the graded algebra of invariants} $(\ss^{\pp}V)^{G}$, the algebra of regular functions on the quotient variety $V/G$. 

The determination of   $(\ss^{\pp}V)^{G}$ goes back to Klein \cite[Ch.\ II]{Kl}. It is generated by 3  homogeneous elements $x,y,z$, subject to one weighted homogeneous relation $F(x,y,z)=0$. Putting  $\mathbf{d}=(\deg x, \deg y, \allowbreak\deg z)$,
we have:

$\bullet$ For $\bar{G}=D_n $,  $\mathbf{d}=(2n+2,2n,4)$,   $F=x^2+y^2z+z^{n+1}$.

$\bullet$ For $\bar{G}=\mathfrak{A}_4$,   $\mathbf{d}=(6,4,4)$,  $F=x^2+y^3+z^3$.

$\bullet$ For $\bar{G}=\mathfrak{S}_4$, $\mathbf{d}=(12,8,6)$, $F= x^2+y^3+z^4$.

$\bullet$ For $\bar{G}=\mathfrak{A}_5$, $\mathbf{d}=(30,20,12)$, $F= x^2+y^3+z^5$.

\smallskip	
\section{Cases with $S(X)=\C$}
\subsection{Varieties with $c_1(X)=0$}
The following result, proved in \cite{Ko}, is a direct consequence of Yau's theorem:
\begin{prop}\label{k}
Let $X$ be a compact K\"ahler variety with $c_1(X)=0$ in $H^2(X,\Q)$, and $\pi _1(X)$ finite. Then $S(X)=\C$.
\end{prop}

With no assumption on $\pi _1(X)$, we know that $X$ is the quotient of a product $A\times Y$, where $A$ is a complex torus and $Y$ is simply connected, by a finite group $G$ acting freely \cite{B2}. It follows that $S(X)$ \emph{is isomorphic to the invariant subring} $(\ss^{\pp}T_0(A))^G$.

\subsection{Varieties of general type}
\begin{prop}\label{gen}
Let $X$ be a  variety of general type. Then $S(X)=\C$.
\end{prop}
This is a consequence of the stronger result that $T_X$ is not pseudo-effective \cite[Proposition 4.11]{HP2}.

\smallskip	
\subsection{Hypersurfaces}
The following result is proved in \cite{HLS}:
\begin{prop}
Let $X\subset\P^{n+1}$ be a smooth hypersurface of degree $d\geq 3$ and dimension $\geq 2$. Then $S(X)=\C$.
\end{prop}
In fact the authors prove the stronger result $H^0(X,\ss^p( T_X(d-3)))=0$, and also that $T_X$ is not pseudoeffective.

\smallskip	
\section{The Krull dimension of S(X)}

A complete description of the ring $S(X)$ is in general intractable, but we  can still  ask for some of its properties,
 for instance its Krull dimension. When $S(X)\neq \C$, it is equal to $1+\kappa (\O_{\P T^*X}(1))$, where $\kappa $ denotes the \emph{Iitaka dimension}  (\cite[Lemma 7.2]{C}). We have $0\leq \dim S(X)\allowbreak \leq 2\dim X$, and all cases can occur.
  In particular,
\[\dim S(X)=2\dim X\ \Longleftrightarrow \ \O_{\P T^*X}(1)\mbox{ big }\Longleftrightarrow\  T_X\mbox{ big}\,.\]
This property holds for toric varieties \cite{Hs}, and also for all rational homogeneous varieties \cite[Corollary 4.4]{GW}. The paper \cite{Li} contains a number of other examples of varieties with a group action whose tangent bundle is big.

\smallskip	
Though the most interesting cases occur when the Kodaira dimension $\kappa (X)$ is $-\infty$, one may ask what can be said when $\kappa (X)\geq 0$. The condition $S(X)\neq \C$, or the weaker condition that $T_X$ is pseudo-effective, imposes strong restrictions on $X$ --- see \cite[Proposition 4.11]{HP2}. The following bound is the main result of this section:
\begin{prop}\label{bound}
$\dim S(X)\leq \dim X-\kappa (X)$. Equality holds if and only if $X$ admits a finite \'etale covering of the form $A\times Y$, where  $A$ is an abelian variety and $Y$ a  variety of general type.
\end{prop}	
It follows in particular that $\dim S(X)>\dim X$ implies $\kappa (X)=-\infty$. 

\smallskip		 
Let us first show that the equality holds when there exists an \'etale covering $A\times Y\rightarrow X$ with $A$ abelian and $Y$ of general type. This follows from Subsection \ref{ab}, Proposition \ref{gen}, and   the following lemma:
\begin{lem}\label{SX}
Let $X,Y$ be smooth projective varieties.

$1)$ We have $S(X\times Y)\cong S(X)\otimes S(Y)$.

$2)$ If $\pi : X\rightarrow Y$ is an \'etale morphism, $\dim S(X)=\dim S(Y)$ and $\kappa (X)=\kappa (Y)$.
\end{lem}
\begin{proof}
1) Let $p_X,p_Y$ be the projections of $X\times Y$ onto $X$ and $Y$. We have $T_{X\times Y}=\allowbreak p_X^*T_X\oplus p_Y^*T_Y$, hence $\ss^{\pp}T_{X\times Y}=p_X^*\ss^{\pp}T_X \otimes  p_Y^*\ss^{\pp}T_Y$. 
The result follows from the K\"unneth formula.

2)   $\pi $ induces a finite \'etale morphism $ T^*X\rightarrow T^*Y$, hence $S(X)=\O(T^*X)$ is a finite algebra over $S(Y)$, thus $\dim S(X)=\dim S(Y)$. The  equality of the Kodaira dimensions is proved in  \cite[Theorem 5.13]{Ue}.
\end{proof}

For the rest of the proof, we will need some preliminary results. 	

\subsection{Slope and positivity of vector bundles}

We fix an ample divisor class $H$ on $X$. We will say that a vector bundle is stable if it is slope-stable with respect to  $H$ ---   same for semi-stability and polystability.  

Let $\mathscr{E}$ be a torsion free coherent sheaf of rank $r$ on $X$. Recall that the \emph{slope} $\mu (\mathscr{E})$ of $\mathscr{E}$ is $\ \dfrac{1}{r} (c_1(\mathscr{E})\cdot H^{n-1})$.  
We denote by $\mu_{\max}(\mathscr{E})$ the maximum of 
$\mu (\mathscr{F})$ for $\mathscr{F}\subseteq \mathscr{E}$, $\mathscr{F}\neq 0$. 

\begin{lem}
	\label{mu}
	Let $E$ and $F$ be two vector bundles on $X$.
	
$1)\ \mu _{\max}(E\otimes F)=\mu _{\max}(E)+\mu _{\max}(F)$.

$2) \ \mu _{\max}(\ss^pE)=p\,\mu _{\max}(E)$.

\noindent In particular, if $E$ and $F$ are semi-stable, then so are $E\otimes F$ and $\ss^q E$ for any $q\geq 1$.
\end{lem}

\begin{proof}
	1) is proved in \cite[Corollary 5.5]{CP}. 
	
	2) Let $\mathscr{F}$ be a subsheaf of $E$ with $\mu (\mathscr{F})=\mu _{\max}(E)$. 
	Then $(\ss^p\mathscr{F})^{**}$ is a subsheaf of $\ss^pE$, hence $\mu _{\max}(\ss^pE)\geq \mu ((\ss^p\mathscr{F})^{**})\geq p\,\mu (\mathscr{F})=p\,\mu _{\max}(E)$. On the other hand since $\ss^pE$ is a subsheaf of $E^{\otimes p}$, we have $\mu _{\max}(\ss^pE)\leq \mu _{\max}(E^{\otimes p})=p\,\mu _{\max}(E)$ by 1), hence 2) holds.
\end{proof}

\subsection{Symmetric algebra of vector bundles}
Let $E$ be a vector bundle of rank $r$ on $X$. We will denote by $S(E)$  the graded algebra $H^0(X,\ss^{\pp} E)$. 
\begin{lem}
	\label{SE}
$1)$ Assume that $E$ is polystable, and $\mu (E)=0$. Then $\dim S(E)\leq r$.

$2)$ Assume $E=F\oplus G$, where $\mu_{\max} (F)\leq 0$ and $\mu _{\max}(G)<0$. Then $S(E)= S(F)$. 
\end{lem}

\begin{proof}
 1) If $E$ is stable and $h^0(E)\neq 0$, there is an injective homomorphism $\O_X\rightarrow E$, which must be an isomorphism; hence $h^0(E)\leq 1$. It follows that $h^0(E)\leq r$ if $E$ is polystable. Now $\ss^qE$ is also polystable \cite[Theorem 3.2.11]{HL}, so $h^0(\ss^qE)\leq \rk \ss^qE=\allowbreak \binom{q+r-1}{r-1} $, hence $\kappa (\O_{\P(E)}(1))\leq r-1$ (see e.g. \cite[Corollary 2.1.38]{La}) and $\dim S(E)\leq r$.

\smallskip	
2) By Lemma \ref{mu} we have, for $p,q\in\N$, $q>0$:
\[ \mu _{\max}(\ss^pF\otimes \ss^q G)= p\, \mu _{\max}(F)+q\,\mu _{\max}(G)< 0\,, \mbox{ hence }\ H^{0}(\ss^pF\otimes \ss^q G)=0\,.\]Therefore $H^0(\ss^pE)=H^0(\ss^pF)$, and $S(E)=S(F)$. 
\end{proof}

\subsection{Proof of Proposition \ref{bound}}
Without loss of generality, we may assume $\kappa(X)\geq 0$ and $\dim S(X)\geq 1$. In particular, the projective manifold $X$ is not uniruled and $T_X$ is pseudo-effective. Moreover, since $\dim S(X)$ and $\kappa(X)$ are invariant under finite \'etale covering (Lemma \ref{SX}), we may replace $X$ by  any finite \'etale covering.

 Proposition 4.11 of \cite{HP2} provides a decomposition 
 \begin{equation}\label{dec}
T_X=F\oplus G
\end{equation}
 where $F$ and $G$ are integrable subbundles, $c_1(F)=0$, and the restriction of $G^*$ to a general curve complete intersection of hypersurfaces in $\abs{mH}$, for $m\gg 0$, is ample. Since a quotient of an ample bundle is ample, this implies $\mu (\mathscr{F})<0$ for any nonzero subsheaf  $\mathscr{F}\subset G$, hence $\mu _{\max}(G)<0$. Then by Lemma \ref{SE} the algebra $S(X)$ is isomorphic to $S(F)$.  
By  \cite[Lemma 2.1]{PT}, $F$ is polystable, hence  Lemma \ref{SE} implies
 \[\dim S(X)=\dim S(F)\leq \rk F \,.\]
By \cite[Proposition 2.6]{PT}, $\det F$
 is a torsion line bundle;
  passing to a finite \'etale covering we may assume $\det F=\O_X$, so that $\det G^*\cong K_X$.  The natural inclusion $G^*\subset \Omega_X^1$ induces an inclusion $\det G^*\subset \Omega_X^k$, where $k=\rk G$. Then the  Bogomolov inequality (\cite[Theorem 4]{Bo}) gives
\[
	\kappa(X)=\kappa(\det G^*)\leq k=\rk G,
	\] hence
\[
	\dim S(X) = \dim S(F) \leq \rk F = \dim X - \rk G \leq \dim X-\kappa(X),
	\]which proves our bound. 
	
Suppose that the equality holds. Then $\dim S(F)=\rk F$ and $\kappa (\det G^*)=k$. By \cite[Lemma 12.4]{Bo}, the latter condition implies that there exists a rational map $f:X\dashrightarrow Y$ to a $k$-dimensional projective manifold such that $\det G^*\subset \Omega_X^k$ coincides with the saturation of the subsheaf $ f^*K_Y\subset \Omega^k_X$. This implies that the foliation $F\subset T_X$ is induced by $f$ and thus  is a regular algebraically integrable foliation.  Since $\det F\cong \O_X$, by the global version of the Reeb stability theorem \cite[Theorem 8.1]{D3}, after replacing $X$ by a finite \'etale covering, we may assume that $X$ is a product $Z\times Y$, with $F= \operatorname{pr} _Z^* T_Z$ and $G\cong \operatorname{pr} _Y^* T_Y$. In particular, we obtain
	\[
	\dim(Y)=\kappa(X,\det G^*)=\kappa(Y)
	\]
 hence $Y$ is of general type.
	
	Finally we use the first condition $\dim S(F)=\rk F$. Since $S(F)$ is canonically isomorphic to $S(Z)$, we get $\dim S(Z)=\dim Z$. Since $c_1(F)=0$, we have $c_1(Z)=0$, hence $Z$ admits a finite \'etale covering of the form $A\times T$, where $A$ is an abelian variety and $T$  a simply connected smooth projective variety with $c_1(T)=0$ \cite{B1}. By Proposition \ref{k} and Lemma \ref{SX}  we have $S(Z)\cong S(A)$, hence $\dim Z=\dim S(Z)=\dim(A)$ (Subsection \ref{ab}), so that $X=Z\times Y$ admits  a finite \'etale covering 
by $A\times Y$.\qed

\smallskip	
\section{Pseudo-effective tangent bundle}

We discuss in this section the structure of non-uniruled projective manifolds $X$ with pseudo-effective tangent bundle. 

\begin{lem}
	\label{l.propertiesPseff}
	$1)$ Let $D$ be a big divisor on $X$. A vector bundle $E$ is pseudo-effective if and only if for any $c>0$, there exist positive integers $i$ and $j$ such that $i>cj$ and 
	\[
	H^0(X,\ss^i E\otimes \O_X(jD))\not=0.
	\]
	
	$2)$ If $E$ is a pseudo-effective vector bundle, then $\mu_ {\max}(E) \geq 0$ for any polarization $H$.
	
	$3)$ Let $F\rightarrow E$ be an injective map of  vector bundles. If $F$ is  pseudo-effective,  $E$  is pseudo-effective.
	
	$4)$ Let $f:Y\rightarrow X$ be a surjective morphism between smooth projective varieties, and let $E$ be a vector bundle on $X$. Then $E$ is pseudo-effective if and only if $f^*E$ is pseudo-effective. 
		
	$5)$ Let $X=Y\times Z$ be a product of smooth projective varieties. Then $T_X$ is pseudo-effective if and only if one of $T_X$ and $T_Z$ is pseudo-effective.
\end{lem}

\begin{proof}
	1) is proved in \cite[Lemma 2.2]{HLS}. 
	
	\smallskip	
	2) If $H^0(X,\ss^i E\otimes \O_X(jD))\not= 0$,   there is an inclusion $\O_X(-jD)\subset \ss^i E$. By Lemma \ref{mu}, we have
	\[
	\mu_ {\max}(E) = \frac{1}{i} \mu_{\max}(\ss^i E) \geq -\frac{1}{i}(jD\cdot H^{n-1}) > -\frac{1}{c} (D\cdot H^{n-1}). 
	\]
	As $c$ is arbitrary, we obtain $\mu_{\max}(E)\geq 0$, which proves 2). 
	
	\smallskip	
	3) follows from 1) and the natural inclusion $\ss^iF\otimes \O_X(jD)\subset \ss^iE\otimes \O_X(jD)$.
	
	\smallskip	
	4) Assume first $\rk E=1$. We only need to show that if $f^*E$ is pseudo-effective, then so is $E$ itself. Indeed, assume  the opposite. By \cite[Theorem 0.2]{BDPP}, there exists a  covering family $\{C_t\}_{t\in T}$ of curves such that $(c_1(E)\cdot C_t)<0$. Let   $\{C_{t'}\}_{t'\in T'}$ be   a covering family of curves on $Y$ such that a general curve $C_{t'}$ is mapped onto some $C_t$. Then we have $(c_1(f^*E)\cdot C_{t'})<0$ by the projection formula, so $f^*E$ is not pseudo-effective by \cite[Theorem 0.2]{BDPP}.
	
	If $\rk E>1$, $f$ induces a surjective morphism $\bar{f}:\P(f^*E^*)\rightarrow \P(E^*)$ such that $\bar{f}^*\O_{\P(E^*)}(1)\cong \O_{\P(f^*E^*)}(1)$; 4) follows from the previous result applied to $\bar{f}$.

	\smallskip	
	5) By 3) and 4), if $T_Y$ or $T_Z$ is pseudo-effective so is $T_X$.  Assume that $T_X$ is pseudo-effective.  Let $H_Y$ and $H_Z$ be ample line bundles on $Y$ and $Z$, respectively.  Then $H:= H_Y\boxtimes H_Z$ is ample. By 1), for any $c>0$, there exist positive integers $i$ and $j$ such that $i>2c j$ and
		\[
	H^0(X,\ss^i T_X\otimes H^j) = H^0(X,\ss^i(T_Z\boxtimes T_Y)\otimes H^j)\not=0.
	\]
	By restricting to   $Y\times \{z\} $ and $\{y\}\times Z $, for $y,z$ general, it follows that there exist non-negative integers $p$ and $q$ such that $p+q=i$, $H^0(Y,\ss^p T_Y\otimes H_Y^j)\not=0$ and \allowbreak$H^0(Z,\ss^q T_Z\otimes H^j_Z)\not=0$. Moreover, as $p+q>2cj$, we also have either $p>cj$ or $q>cj$. Since $c$ is arbitrary and $H$ is ample, it follows from 1) that one of $T_Z$ and $T_Y$ is pseudo-effective.
\end{proof}

\noindent\emph{Remark}$.-$ In general, if the tangent bundle $T_X$ of a smooth projective variety $X$ is pseudo-effective and splits into a direct sum $F\oplus G$ of vector bundles, it is  not clear to us  whether one of $F$ or $G$ is pseudo-effective. Indeed, the splitting of $T_X$ in general does not imply the splitting of $X$ itself, as simple abelian varieties or Hilbert modular varieties show. However, it is conjectured by the first author in \cite{B3} that this splitting should come from a splitting of the universal cover of $X$.

\bigskip	
Recall that a rank $r$ vector bundle $E$ on $X$ is called \emph{unitary flat} if it is associated to an irreducible representation $\pi _1(X)\rightarrow \operatorname{U}(r) $.

\begin{conj}
	\label{c.PseffTangent}
	Let $X$ be a non-uniruled projective manifold. Then $T_X$ is pseudo-effective if and only if there exists a finite \'etale covering $X'\rightarrow X$ such that $T_{X'}$ contains a nonzero unitary flat subbundle.
\end{conj}

    \noindent\emph{Remarks}$.-$ 1) A unitary flat vector bundle is nef, hence pseudo-effective. So if $T_{X'}$ contains a nonzero unitary flat subbundle, it is pseudo-effective (Lemma \ref{l.propertiesPseff}, 3)), hence $T_X$  is pseudo-effective (Lemma \ref{l.propertiesPseff}, 4)).
    
   \smallskip	 
    2) If the tangent bundle $T_X$ of a non-uniruled projective $X$ contains a unitary flat subbundle $F$, then $F$ is actually a regular foliation with $\det(F)$ torsion by \cite[Lemma 2.1 and Proposition 2.6]{PT}. We refer the reader to \cite{PT} for more discussion on the structure of this kind of foliations.
		
    \smallskip	
    3) Very recently J. Jia, Y. Lee and G. Zhong have studied in \cite{JLZ} the non-uniruled smooth projective surfaces $S$ with pseudo-effective tangent bundle. They prove that up to a finite \'etale covering, $S$ is either an abelian surface or a product $E\times C$ of an elliptic curve $E$ and a curve $C$ of genus $\geq 2$. This solves Conjecture \ref{c.PseffTangent} in dimension two. 
		
	In higher dimension, it is asked in \cite[Question 1.2]{JLZ} whether the pseudo-effectivity of the tangent bundle of an $n$-dimensional non-uniruled projective manifold $X$ is equivalent to $c_n(X)=0$ and $\widehat{q}(X)>0$, where $\widehat{q}(X)$ is the \emph{augmented irregularity} of $X$. The answer is negative in general. For instance, let $X=Y\times Z$ be the product of an irreducible simply connected Calabi-Yau variety $Y$ with vanishing top Chern class\footnote{See for instance \cite[p. 1221]{KS} for the construction of  threefolds with this property.} and a variety $Z$ of general type with $q(Z)>0$. The tangent bundles of $Y$ and $Z$ are not pseudo-effective (see \cite[Theorem 1.6]{HP1} and \cite[Proposition 4.11]{HP2}). So  Lemma \ref{l.propertiesPseff} says that $T_X$ itself is not pseudo-effective.
	
\bigskip

Because of  the decomposition \eqref{dec},  Conjecture \ref{c.PseffTangent} is closely related to the following conjecture proposed by J.V. Pereira and F. Touzet in \cite[\S\,6.5]{PT}:

\begin{conj}
	\label{c.PTConj}
	Let $X$ be a non-uniruled projective manifold, and let $F\subsetneq T_X$ be a regular foliation such that $F$ is stable for some polarization, $c_1(F)=0$, and $c_2(F)\neq 0$. Then $F$ is algebraically integrable.
\end{conj}

\begin{prop}
	Assume that  Conjecture \emph{2} holds for $\dim(X)\leq n$. Then  Conjecture \emph{1} holds for $\dim(X)\leq n$.
\end{prop}

\begin{proof}
	Assume that $T_X$ is pseudo-effective. Let $T_X=F\oplus G$ be the decomposition  \eqref{dec}. Then $F$ is a regular foliation with $c_1(F)=0$. By \cite[Theorem 6.9]{D1}, there exist complex projective manifolds $Y$ and $Z$, a finite \'etale cover $\pi :Y\times Z \rightarrow X$, and a regular foliation $H$ on $Y$ with $c_1(H)=c_2(H)=0$ such that $\pi ^*F=p_Y^*H\oplus p_Z^*T_Z$. Since $H$ is polystable  \cite[Lemma 2.1]{PT}, it is a direct sum of  unitary flat bundles \cite[Corollary 8.1]{UY}. Therefore it suffices to prove that $H\neq 0$.

Since $c_1(F)=0$ and $c_1(H)=0$, we get $c_1(Z)=0$. Therefore
 there exists a finite \'etale covering $A\times T\rightarrow Z$,  where $A$ is an abelian variety and $T$ is a simply connected smooth projective variety with $c_1(T)=0$  \cite{B1}. Without loss of generality, we may assume that $Z=A\times T$. Moreover, after replacing $Y$ by $A\times Y$, we may assume in addition that $Z$   is simply connected. In particular, the tangent bundle $T_Z$ is not pseudo-effective  \cite[Theorem 1.6]{HP1},  so $T_Y$ is pseudo-effective by Lemma \ref{l.propertiesPseff}. 
	
	Applying \cite[Theorem 2.2]{PT} to $Y$ yields   a regular foliation $J$ on $Y$ such that $T_Y=H\oplus J$. We have on one hand $\pi ^*T_X=\pi ^*F\oplus \pi ^*G$, and on the other hand \[\pi ^*T_X=p_Y^*(H\oplus J)\oplus p_Z^*T_Z=\pi ^*F\oplus p_Y^*J\,;\] this implies $p_Y^*J\cong \pi ^*G$. Since $\mu _{\max}(G)<0$, $J$ is not pseudo-effective (Lemma \ref{l.propertiesPseff}). Therefore $H\neq 0$
	 and we are done.
\end{proof}

Conjecture \ref{c.PTConj} is  wide open in general.  It is known in  the following    cases, proved by F. Touzet and S. Druel (\cite{To} and \cite{D1}).

\begin{prop}
	Conjecture \emph{2} holds if $\rk(F)\leq 3$ or $\rk(F)=\dim(X)-1$. In particular, it holds for $\dim(X)\leq 5$.
\end{prop}

\begin{proof}
	If $\rk(F)\leq 3$, this is proved in \cite[Proposition 6.8]{D1}. Assume  $\rk(F)=\dim(X)-1$, and that  $F$ is not algebraically integrable.
	By \cite[Th\'eor\`eme 1.2]{To},   there exists an abelian variety $A$, a smooth projective variety $Y$ with $c_1(Y)=0$, a finite \'etale covering $\pi :A\times Y\rightarrow X$ and a linear foliation $H$ on $A$ such that $\pi ^*F=p_A^*H\oplus p_Y^*T_Y$. Since $F$ is stable for some polarization,  Proposition 8.1 of \cite{D2} implies that $Y$ is a point.  Then $\pi ^*F=H$ is trivial,  hence $c_2(F)=0$, a contradiction.
\end{proof}

\begin{cor}
	Conjecture \emph{1} holds for $\dim(X)\leq 5$.
\end{cor}

\end{document}